\documentclass[12pt, a4paper,twosided,reqno]{amsart}
\usepackage{enumerate, cite}
\usepackage{hyperref}
\usepackage{amsmath}

\newtheorem{theorem}{Theorem}
\newtheorem{lemma}{Lemma}

\newtheorem{example}{Example}
\newtheorem{corollary}{Corollary}

\newtheorem*{thA}{Theorem A}
\newtheorem*{thB}{Theorem B}
\newtheorem*{thC}{Theorem C}
\newtheorem*{thD}{Theorem D}

\allowdisplaybreaks
\author[G. Pant, M. Saini]{Garima Pant and Manisha Saini}
\address{Garima Pant; department of mathematics, university of delhi, delhi-110007, india.}
\email{garimapant.m@gmail.com}
\address{department of mathematics, university of delhi, delhi-110007, india.}
\email{msaini@maths.du.ac.in, sainimanisha210@gmail.com}
\thanks {Research work of the first author is supported by research fellowship from University Grants Commission (UGC), New Delhi, India.}
\thanks {The second author is Senior Research Fellow (UGC, New Delhi, India).}
\title[study of solutions...]{ Study of Solutions of Certain kind of Non-Linear Differential Equations}
\begin{document}
\maketitle
\begin{abstract}
In this paper we study about the existence of solutions of certain kind of non-linear differential and differential-difference equations. We give partial answer to a problem which was asked by chen et al. in \cite{chw}.\\
\end{abstract}

 We assume that the readers are familiar with the standard notations of Nevanlinna theory such as characteristic function $T(r,f)$, proximity function $m(r,f)$, integrated counting function $N(r,f)$, unintegrated counting function $n(r,f)$ and first main theorem for a meromorphic function $f$ which is defined in the complex plane \cite{ilaine, hayman, yanglo}. To make our paper self contained, we present the elementary definitions of the order of growth $\rho(f)$ and exponent of convergence of zeros $\lambda(f)$ of a meromorphic function $f$
$$\rho(f)=\limsup_{r\to\infty}\frac{\log^{+} T(r,f)}{\log r}$$
 and
 $$\lambda(f)=\limsup_{r\to\infty}\frac{\log^{+} n(r,1/f)}{\log r},$$
respectively. If $f$ is a meromorphic function, then $S(r,f)$ denotes such quantities which are of growth $o(T(r,f))$ as $r\to\infty$, outside of a possible exceptional set of finite linear measure. A meromorphic function $h(z)$ is said to a small function of $f(z)$ if $T(r,h)=S(r,f)$ and vice versa. Note that finite sum of the quantities which are of $S(r,f)$ kind is again $S(r,f)$. A differential polynomial in $f$ is a polynomial in $f$ and its derivatives with small function of $f$ as its coefficients. For a given meromorphic function $f$ and a constant $c$, $f(z+c)$ is known as a shift of $f$.\\
To study solubility and existence of solutions of non-linear differential equations or differential-difference equations are interesting but also difficult. Many researchers have been studied these equations, one may refer to see \cite{li,li1,ly,rx,ww}. \\
In this paper, our aim is to present some studies on the solutions of certain kind of non-linear differential equations and differential-difference equations which have been studied even before. 
In 2021, Chen et al.\cite{chw} studied the following non-linear differential-difference equation
\begin{equation}\label{chweq}
f^n(z)+wf^{n-1}(z)f'(z)+q(z)e^{Q(z)}f(z+c)=u(z)e^{\nu(z)},
\end{equation}
where $n$ is a natural number, $w$ and $c\neq 0$ are constants. $q, Q, u$ and $\nu$ are  non constant polynomials.\\

In our first result, we study the same non-linear differential-difference equation \eqref{chweq} and prove the following result:
\begin{thA}\label{mainth1}	
Suppose that $n\geq 3$ for $w\neq 0$ and $n\geq 2$ for $w=0$. Let $f$ be a finite order transcendental entire solution of equation \eqref{chweq} satisfying $\lambda(f)<\rho(f)$. Then $f$ satisfies one of the following:
\begin{enumerate}[$(i)$]
\item $\rho(f)<\deg\nu=\deg Q$ and $f=Ce^{-z/w}$, where $C$ is a non-zero constant.
\item $\rho(f)=\deg\nu=\deg Q$ 
\end{enumerate}
\end{thA}
Note that every non-vanishing transcendental entire function $f$ of non-zero finite order satisfies $\lambda(f)<\rho(f)$. Therefore we get an immediate corollary from the above theorem.
\begin{corollary}\label{partialans1}
Suppose that $n\geq 3$ for $w\neq 0$ and $n\geq 2$ for $w=0$. If $f$ is a non-vanishing transcendental solution of non-zero finite order of equation \eqref{chweq}. Then
$f$ satisfies the conclusion of Theorem A. 
\end{corollary}
Prior to Theorem A \rm\ref{mainth1}, Chen et al.\cite{chw} proved the following result:

\begin{theorem}
Assume $n\geq 3$ for $w\neq 0$ and $n\geq 2$ for $w=0$. Let $f$ be a non-vanishing transcendental solution of finite order of equation \eqref{chweq}.	Then each 
$f$ satisfies one of the following:  
\begin{enumerate}[$(i)$]
\item $\rho(f)<\deg\nu=\deg Q$ and $f=Ce^{-z/w}$, where $C$ is a constant.
\item $\rho(f)=\deg Q\geq \deg \nu$ 
\end{enumerate} 
\end{theorem}

After stating the above theorem, Chen et al.\cite{chw}   posed a problem namely Problem $1.10$,  in which they 
asked ``Can the conclusion $(ii)$ in the above theorem be further improved to $\rho(f)=\deg Q= \deg \nu$?"\\
In the Corollary \ref{partialans1}, we have given partial answer to that problem. \\

Chen et al.\cite{chw} also studied a non-linear differential-difference
equation in which left hand side of equation \eqref{chweq} is unaltered and right hand side of the same equation is replaced by $P_{1}e^{\lambda z}+P_{2}e^{-\lambda z}$. Then they provided the following result:
\begin{theorem}
Suppose $f$ is a transcendental entire solution with finite order of equation
\begin{equation}\label{dedegn}
f^{n}(z)+wf^{n-1}(z)f^{'}(z)+q(z)e^{Q(z)}f(z+c)=P_{1}e^{\lambda z}+P_{2}e^{-\lambda z},
\end{equation}
where $n$ is a natural number, $w$ is a constant and $c,P_{1},P_{2},\lambda$ are non-zero constants, $q\not\equiv 0$ is a polynomial and $Q$ is a non-constant polynomial. Then the following conclusions hold:
\begin{enumerate}[(i)]
\item If $n\geq 4$ for $w\neq 0$ and $n\geq 3$ for $w=0$, then each solution $f$ satisfies $\rho(f)=\deg Q=1$.
\item If $n\geq 1$ and $f$ is a solution which belongs to $\Gamma_{0}=\{e^{\alpha(z)}:\alpha(z)$ is a non-constant polynomial\},  then
$$f(z)=e^{(\lambda z/n)+a}, \qquad Q(z)=-\frac{n+1}{n}\lambda z+b$$
or $$f(z)=e^{(-\lambda z/n)+a}, \qquad Q(z)=\frac{n+1}{n}\lambda z+b,$$
where $a$ and $b$ are constants.
\end{enumerate}
\end{theorem}
Now it is natural to ask what can we say about the solutions of equation \eqref{dedegn} when $n=3$ and $w\neq 0$? In this sequence, we prove the following result:

\begin{thB}\label{mainth3}
Suppose that $w,c,P_{1},P_{2}$ and $\lambda$ are non-zero constants, $q(z)$ is a non-zero polynomial, $Q(z)$ is a non constant polynomial. Then there does not exist any finite order transcendental entire function satisfying equation
\begin{equation}\label{dedeg3}
f^{3}(z)+wf^{2}(z)f'(z)+q(z)e^{Q(z)}f(z+c)=P_{1}e^{\lambda z}+P_{2}e^{-\lambda z}.
\end{equation}
with a finite non-zero Borel exceptional value.
\end{thB}
The following example shows that if $f$ is a transcendental entire function satisfying equation \eqref{dedeg3}, then $f$ may have zero as a Borel exceptional value.
\begin{example}
The function $f(z)=2e^{3z}$ satisfies differential-difference equation 
\begin{equation*}
f^{3}+f^{2}f^{'}+e^{12z}f(z+\log 2)=32e^{9z}+16e^{-9z}.
\end{equation*}
Here $f(z)$ has zero as a Borel exceptional value. 
\end{example}
In 2020, Xue \cite{xue} and Chen et al. \cite{cl} studied a new kind of non-linear differential equations and proved the following results:
\begin{theorem}\cite{xue} \label{xueth}
Suppose $n\geq 2$ and $k$ are integers, $P_{n-1}(z,f)$ is an algebraic differential polynomial in $f(z)$ of degree at most $n-1$. Let $P_{i}$ and $\alpha_{i}$ be non-zero constants for $i=1,2,3$, and $|\alpha_{1}|>|\alpha_{2}|>|\alpha_{3}|$. If $f$ is a transcendental entire function  satisfying the following differential equation
\begin{equation}
f^{n}(z)+P_{n-1}(z,f)=P_{1}e^{\alpha_{1}z}+P_{2}e^{\alpha_{2}z}+P_{3}e^{\alpha_{3}z},
\end{equation}
then $f(z)=Ce^{\alpha_{1}z/n}$, where $C$ is a non-zero constant such that $C^{n}=P_{1}$, and $\alpha_{1},\alpha_{2}, \alpha_{3}$ are in one line.
\end{theorem}
Note that an algebraic differential polynomial $P(z,f)$ is a polynomial in $f$ whenever $f$ is a polynomial in $f(z)$ and its derivatives with polynomials as its coefficients.
\begin{theorem}\cite{cl}\label{clth}
Suppose $n\geq 5$ is an integer and $Q_{d}(z,f)$ is a differential polynomial in $f$ of degree $d\leq n-4$ with rational functions as its coefficients. Also suppose $p_{1}(z), p_{2}(z), p_{3}(z)$ are non-zero rational functions and $\alpha_{1},\alpha_{2},\alpha_{3}$ are non constant polynomials such that $\alpha_{1}^{'},\alpha_{2}^{'},\alpha_{3}^{'}$ are distinct to each other. If $f$ is a meromorphic solution with finitely many poles of the following differential equation 
\begin{equation}
f^{n}(z)+Q_{d}(z,f)=p_{1}(z)e^{\alpha_{1}(z)}+p_{2}(z)e^{\alpha_{2}(z)}+p_{3}(z)e^{\alpha_{3}(z)},
\end{equation}
then $\frac{\alpha_{1}^{'}}{\alpha_{2}^{'}}$, $\frac{\alpha_{2}^{'}}{\alpha_{3}^{'}}$ are rational numbers, and $f(z)$ must be of 
$$f(z)=r(z)e^{P(z)},$$
where $r(z)$ is a non-zero rational function and $P(z)$ is a non constant polynomial. Furthermore, there must exist positive integers $m_{1},m_{2},m_{3}$ with $\{m_{1},m_{2},m_{3}\}=\{1,2,3\}$ and distinct integers $k_{1},k_{2}$ with $1\leq k_{1},k_{2}\leq d$ such that $\alpha_{m_{1}}^{'}:\alpha_{m_{2}}^{'}:\alpha_{m_{3}}^{'}=n:k_{1}:k_{2}$, $nP'=\alpha_{m_{1}}^{'}$, $Q_{d}(z,f)\equiv p_{m_{2}}(z)e^{\alpha_{m_{2}}(z)}+p_{m_{3}}(z)e^{\alpha_{m_{3}}(z)}$.
\end{theorem}
After stating the above result, Chen et al. \cite{cl} mentioned that ``in Theorem \ref{clth}, the condition $n\geq 5$ is necessary" in a remark namely 'Remark $1.2$'. By some examples, they have also shown that if $n\leq 4,$ then Theorem \ref{clth} may not hold.
In this sequence we prove the following results:
\begin{thC}\label{mainth2}
Let $n\geq 2$, and $P(z,f)$ be a linear differential polynomial in $f(z)$, its derivatives with coefficients being rational functions. Suppose $P_{j}(z)$ and $\alpha_{j}(z)$ $(j=1,2,3)$ are non-zero rational functions and non constant polynomials, respectively. Also $\alpha_{1}^{'}(z),\alpha_{2}^{'}(z) $and $\alpha_{3}^{'}(z)$ are distinct to each other. If $f$ is a transcendental meromorphic solution of the following differential equation
\begin{equation}\label{dedeg2}
f^{n}(z)+P(z,f)=P_{1}(z)e^{\alpha_{1}(z)}+P_{2}(z)e^{\alpha_{2}(z)}+P_{3}(z)e^{\alpha_{3}(z)},
\end{equation}
	
such that $f$ has finitely many poles, then $f$ is of finite order but not in the form of $f(z)=s(z)e^{q(z)}$, where $s(z)$ is a non-zero rational function and $q(z)$ is a non constant polynomial.
\end{thC}

Note that Theorem \ref{clth} is an improvement of Theorem \ref{xueth} which can be seen by the following examples.  
\begin{example}
The differential equation 
$$f^{2}+f^{''}-2f^{'}-f=2e^{2z}+2e^{3z}+e^{4z}$$
has a solution $f(z)=1+e^{z}+e^{2z}$ which is not in the required form as in Theorem $\ref{xueth}$.
\end{example}
\begin{example}
The differential equation 
$$f^{3}-2f^{'}-4f=e^{9z}+2e^{3z}+6e^{z}$$
has a solution $f(z)=e^{3z}+2$ which is not in the required form as in Theorem $\ref{xueth}$.
\end{example}
\begin{example}
The differential equation
$$f^{4}+\frac{3}{2}f^{'}+f=16e^{8z}-32e^{6z}+23e^{4z}$$ 
has a solution $f(z)=2e^{2z}-1$ which is not in the required form as in Theorem $\ref{xueth}$.
\end{example}
We observe that the above examples also hold Theorem C. Next, we combine Theorem \ref{clth}, Theorem C and we have an immediate result:
\begin{corollary}
Let $n\geq 5$, and $P(z,f)$ be a linear differential polynomial in $f(z)$, its derivatives with rational functions as its coefficients. Suppose $P_{j}(z)$ and $\alpha_{j}(z)$ $(j=1,2,3)$ are non-zero rational functions and non constant polynomials, respectively. Also $\alpha_{1}^{'}(z),\alpha_{2}^{'}(z),\alpha_{3}^{'}(z)$ are distinct to each other. Then there exist no meromorphic solution having finitely many poles which satisfies equation \eqref{dedeg2}.
\end{corollary}
Before stating next result we first define that $Q_{2}^{\star}(z,f)$ is a differential polynomial of degree $2$ with rational functions $a_{i}(z)$ as its coefficients and at least one coefficient of $f^{(l)} (l\geq 0)$ must be non-zero.
\begin{thD}
Suppose $P_{1},P_{2},P_{3}$ are non-zero constants, and $\alpha_{1},\alpha_{2},\alpha_{3}$ are non-zero and distinct constants. If $f$ is an entire solution of the following differential equation
\begin{equation}\label{dedeg3sp}
f^{n}(z)+Q_{2}^{\star}(z,f)=P_{1}e^{\alpha_{1}z}+P_{2}e^{\alpha_{2}z}+P_{3}e^{\alpha_{3}z},
\end{equation}
where $n\geq 3$ and satisfying $\lambda(f)<\rho(f)$. Then $\rho(f)=1$ and $f(z)$ must be in the form of $f(z)=Ce^{az}$, where $C$ and $a$ are non-zero constant. Particularly, we can say that there exist positive integers $k_{1},k_{2},k_{3}$ with $\{k_{1},k_{2},k_{3}\}=\{1,2,3\}$ such that $\alpha_{k_{1}}:\alpha_{k_{2}}:\alpha_{k_{3}}=n:2:1$, $na=\alpha_{k_{1}}$ and $C=P_{k_{1}}^{1/n}$.
\end{thD}
By the following examples, we show that our hypothesis keeps necessary conditions to conclude the result.
\begin{example}
The differential equation
$$f^{3}+3f^{2}=e^{3z}-6e^{2z}+9e^{z}$$
has a solution $f(z)=e^{z}-3$, which is not in the required form, here $3f^{2}$ is not in the form of $Q_{2}^{\star}(z,f)$.
\end{example}
\begin{example}
The differential equation
$$f^{3}+2f^{2}+f^{'^{2}}=e^{3z}-3e^{2z}+4e^{z}$$
has a solution $f(z)=e^{z}-2$, which is not in the required form, here $2f^{2}+f^{'^{2}}$ is not in the form of $Q_{2}^{\star}(z,f)$.
\end{example}
\begin{example}
The differential equation
$$f^{3}+ff^{'}-f=e^{3z}-2e^{2z}+e^{z}$$
has a solution $f(z)=e^{z}-1$, which is not in the required form and $\lambda(f)=\rho(f)$.
\end{example}
\begin{example}
The differential equation
$$f^{4}-f^{2}-5f^{'^2}=e^{4z}-4e^{3z}-2e^{z}$$
has a solution $f(z)=e^{z}-1$, which is not in the required form and given coefficients of exponent in the right hand side of the above equation don't preserve required ratio, here $-f^{2}-5f^{'^2}$ is not in the form of $Q_{2}^{\star}(z,f)$.
\end{example}
Note that Theorem D is an improvement of Theorem \ref{xueth}.
\section{\textbf{Preliminary Results}}
In this section, we collect those lemmas which we use to prove our desired results.\\

The following lemma is related to the proximity function of logarithmic derivative of a meromorphic function $f$. 

\begin{lemma}\cite{ilaine} \label{il}
Suppose $f$ is a transcendental meromorphic function and $k\geq1$ is an integer. Then
$$m\left(r,\frac{f^{(k)}}{f}\right)=S(r,f).$$
\end{lemma}
The following lemma estimates the characteristic function of a shift of a meromorphic function $f$.
\begin{lemma}\cite{cf}\label{cflemma}
Suppose $f$ is a meromorphic function of finite order $\rho$ and $c$ is a non-zero complex constant. Then for every $\epsilon$,
$$T(r, f(z+c))=T(r,f)+O(r^{\rho-1+\epsilon})+O(\log r).$$
\end{lemma}
Next lemma is the difference analogue of the lemma on the logarithmic derivative of a meromorphic function $f$ having finite order.
\begin{lemma}\cite{cf}\label{hk}
Suppose $f$ is a meromorphic function with $\rho(f)<\infty$ and $c_1,c_2\in\mathbb{C}$ such that $c_1\neq c_2$, then  for each $\epsilon>0$, we have
$$m\left(r,\frac{f(z+c_1)}{f(z+c_2)}\right) =O(r^{\rho-1+\epsilon}).$$
\end{lemma}
The following lemma estimates the characteristic function of an exponential polynomial $f$. This lemma can be seen in \cite{whl}.
\begin{lemma}\label{whllem}
Suppose $f$ is an entire function given by
$$f(z)=A_{0}(z)+A_{1}(z)e^{w_{1}z^{s}}+A_{2}(z)e^{w_{2}z^{s}}+...+A_{m}(z)e^{w_{m}z^{s}},$$
where $A_{i}(z);0\leq i\leq m$ denote either exponential polynomial of degree $<s$ or polynomial in $z$, $w_{i};1\leq i\leq m$ denote the constants and $s$ denotes a natural number. Then
$$T(r,f)=C(Co(W_{0}))\frac{r^{s}}{2\pi}+o(r^{s}),$$
Here $C(Co(W_{0}))$ is the perimeter of the convex hull of the set $W_{0}=\{0,\tilde{w_{1}},\tilde{w_{2}},...,\tilde{w_{m}}\}$.
\end{lemma}


The following lemma plays a key role to prove all our results in this article.
\begin{lemma}\label{imp}\rm\cite{yybook}
Suppose $f_1,f_2,...,f_n (n\geq2)$ are meromorphic functions and $h_1,h_2,...,h_n$ are entire functions satisfying 
\begin{enumerate}
\item $\sum_{i=1}^{n}f_ie^{h_i}\equiv 0$.
\item For $1\leq j<k\leq n$, $h_j-h_k$ are not constants .
\item For $1\leq i\leq n, 1\leq m<k\leq n$,\\
$T(r, f_i)=o(T(r,e^{(h_m-h_k)}))$ as $r\to\infty$, outside a set of finite logarithmic measure.
\end{enumerate}
Then $f_i\equiv 0$ $(i=1,2,...,n).$
	
\end{lemma}
\section{\textbf{Proof of Theorems}}

\begin{proof}[\textbf{\underline{Proof of Theorem A}}]
 Given that $f$ is a transcendental entire solution of finite order of equation \eqref{chweq} satisfying $\lambda(f)<\rho(f)$. \\
 \textbf{Case 1:} First, we consider $n\geq 3$ for $w\neq 0$, then using  Weierstrass factorisation theorem, $f$ can be written as
\begin{equation}\label{fbywft}
f(z)=\alpha(z)e^{\beta(z)},
\end{equation}
where $\beta(z)$ is a non constant polynomial of degree $t$ and $\alpha(z)$ is an entire function satisfies $\rho(\alpha)=\lambda(f)<\rho(f)=t$.\\
For simplicity, rewriting \eqref{chweq} as 
\begin{equation}\label{chweqsimple}
f^n+wf^{n-1}f'+qe^{Q}f_c=ue^{\nu},
\end{equation}
where $f_c=f(z+c)$.
Now, we study the following three sub cases.\\
\textbf{Case 1.1:}
When $\rho(f)<\deg Q$, then we have $T(r,f)=S(r,e^{Q})$.\\
 Also, it is well known that a meromorphic function and its derivative have same order, thus we have $T(r,f')=S(r,e^{Q})=T(r,f)$. Using Lemma \ref{cflemma}, we obtain $$T(r,f_c)=S(r,e^{Q})$$. \\
From equation \eqref{chweqsimple}, we get
\begin{align*}
T(r,ue^{\nu})&=T(r,f^n+wf^{n-1}f'+qe^{Q}f_c)\\
&\leq T(r,f^{n-1})+T(r,f+wf')+T(r,q)+T(r,e^{Q})+T(r,f_c)+\log 2\\
&\leq nT(r,f)+T(r,f')+T(r,e^Q)+T(r,f_c)+T(r,q)+O(1)\\
&= T(r,e^Q)+S(r,e^Q) \qquad \mbox{(Since $T(r,f)=T(r,f')=T(r,f_c)=S(r,e^Q)$)}
\end{align*} 
and
\begin{align*}
T(r,e^Q)&=T\left(r,\frac{qf_ce^{Q}}{qf_c}\right)\\
&\leq T(r, qf_ce^{Q}+wf^{n-1}f'+f^{n})+T(r,wf^{n-1}f'+f^{n})+T\left(r,\frac{1}{qf_c}\right)+O(1)\\
&\leq T(r,ue^{\nu})+nT(r,f)+T(r,f')+T(r,f_c)+T(r,q)+O(1)\\
&=T(r,ue^{\nu})+S(r,e^{Q}).
\end{align*}
From the last two inequalities, 
$$T(r,e^Q)=T(r,ue^{\nu})+S(r,e^{Q}).$$
This implies $\deg Q=\deg \nu$. Now Differentiation of equation \eqref{chweqsimple} gives
\begin{equation}\label{maineq2}
nf^{n-1}f'+w(n-1)f^{n-2}(f')^2+wf^{n-1}f^{''}+Ae^{Q}=(u'+u\nu')e^{\nu},
\end{equation}
where, $A=q'f_c+qQ'f_c+qf'_c$. Eliminating $e^{\nu}$ from equation \eqref{chweqsimple} and \eqref{maineq2}, then we obtain
$$D_1e^{Q}+D_2=0,$$
where $$D_1=(u'+u\nu')qf_c-uA,$$
 $$D_2=(u'+u\nu')(f^n+wf^{n-1}f')-u(nf^{n-1}f'+w(n-1)f^{n-2}(f')^2+wf^{n-1}f").$$
Using Lemma \ref{imp}, we have $D_1=D_2\equiv 0$.\\
From $D_1=0$, we get 
$$\frac{u'}{u}+\nu'=\frac{q'}{q}+\frac{f'_c}{f_c}+Q'.$$
On doing integration of the above equation, we have
$$qf_ce^Q=C_1ue^{\nu}.$$
where $C_1$ is a non-zero constant. If $C_1=1$, then $qf_ce^Q=ue^{\nu}$ and by equation \eqref{chweqsimple}, we obtain
$$f^{n}+wf^{n-1}f'=0.$$
This yields that $f=Ce^{-z/w}$, where $C$ is a non-zero constant. This is the required conclusion $(i)$.\\
If $C_1\neq 1$, then $qf_ce^Q=C_1ue^{\nu}$.
This implies that $f=\frac{C_1u_{-c}e^{\nu_{-c}-Q_{-c}}}{q_{-c}},$
where $u_{-c}=u(z-c)$, $\nu_{-c}=\nu(z-c)$ and $Q_{-c}=Q(z-c)$.\\
Substituting the above value of $f$ into the equation \eqref{chweqsimple}, we have
\begin{equation*}
\begin{split}
\left(\frac{C_1u_{-c}}{q_{-c}}\right)^{n-1}\left[\frac{C_1u_{-c}}{q_{-c}}+C_1w\left(\left(\frac{u_{-c}}{q_{-c}}\right)^{'}+\left(\frac{u_{-c}}{q_{-c}}(\nu_{-c}-Q_{-c})^{'}\right) \right)\right]&e^{n(\nu_{-c}-Q_{-c})}=\\
&(1-C_1)ue^{\nu}.
\end{split}
\end{equation*}
Since $\deg Q=\deg \nu>\rho(f)=\deg \beta\geq 1$ and $\rho(f)=\deg(\nu_{-c}-Q_{-c})$. Using Lemma \ref{imp}, we have $(1-C_1)u\equiv 0$. This implies either $C_1=1$ or $u=0$, which is not possible.\\
\textbf{Case 1.2:} When $\rho(f)=\deg Q$, we substitute the value of $f$ from equation \eqref{fbywft} into equation \eqref{chweqsimple}
\begin{equation}\label{maineq1}
(\alpha^{n}+w\alpha^{n-1}(\alpha^{'}+\beta^{'}\alpha))e^{n\beta}+q\alpha_ce^{Q+\beta_c}=ue^{\nu},
\end{equation}
where $\beta_c=\beta(z+c)$. Now, we consider
$\beta(z)=a_tz^{t}+a_{t-1}z^{t-1}+...+a_0, (a_t\neq 0)$ and $Q(z)=b_tz^{t}+b_{t-1}z^{t-1}+...+b_0,(b_t\neq 0)$ are polynomials of degree $t$ provided $\rho(f)=\deg \beta=t$.\\
Let $\deg \nu<\deg Q=\rho(f)=\deg \beta=t$, then from equation \eqref{maineq1}, we have
\begin{equation}\label{mainsubeq1}
[\alpha^{n}+w\alpha^{n-1}(\alpha'+\alpha \beta')]e^{\alpha_1}e^{na_tz^{t}}+q\alpha_ce^{\alpha_2}e^{(b_t+a_t)z^{t}}=ue^{\nu},
\end{equation}
where $\alpha_1=n(a_{t-1}z^{t-1}+...+a_0)$ and $\alpha_2=(b_{t-1}+cta_{t}+a_{t-1})z^{t-1}+...+(b_{0}+ca_{t}+a_{0})$ are polynomials of degree at most $t-1$.\\
If $a_t\neq \pm b_t$ , then using Lemma \ref{imp}, we get $u\equiv 0$, which is not possible.\\
If $a_t=b_t$, then again using Lemma \ref{imp}, we get the same contradiction.\\
If $a_t=-b_t$, then
$$[\alpha^{n}+w\alpha^{n-1}(\alpha'+\alpha \beta')]e^{\alpha_1}e^{na_tz^{t}}=ue^{\nu}-q\alpha_ce^{\alpha_2}.$$
Using Lemma \ref{imp}, we get 
$$\alpha^{n-1}(\alpha+w(\alpha'+\alpha \beta'))\equiv 0.$$
This gives $\alpha=C_2e^{-z/w-\beta}$, where $C_2$ is a non-zero constant, which implies $\rho(\alpha)=\deg \beta=\rho(f)$. This is  a contradiction to the fact that $\rho(\alpha)<\rho(f)$. Thus $\deg \nu<\deg Q$ is not possible.\\
Let $\deg \nu>\deg Q=\rho(f)=t$ and we consider $\nu(z)=\nu_{s}z^{s}+\nu_{s-1}z^{s-1}+...+\nu_{0}$ such that $s>t$, then from equation \eqref{maineq1}, we have
$$[\alpha^{n}+w\alpha^{n-1}(\alpha'+\alpha \beta')]e^{\alpha_1}e^{na_tz^{t}}+q\alpha_ce^{\alpha_2}e^{(b_t+a_t)z^{t}}=ue^{\alpha_3}e^{\nu_{s}z^{s}},$$
where $\alpha_1$, $\alpha_2$ are same polynomials as in equation \eqref{mainsubeq1} and $\alpha_3=\nu_{s-1}z^{s-1}+...+\nu_{0}$ is a polynomial of degree at most $s-1$.\\
If $a_t\neq \pm b_t$, then using Lemma \ref{imp}, we get $u\equiv 0$, which is not possible.\\
 If $a_t=b_t$, then again using Lemma \ref{imp}, we get  $u\equiv 0$, which is not possible.\\
 If $a_t= -b_t$, then we also get same contradiction as for the above cases. Therefore $\deg \nu>\deg Q$ is not possible.
Hence $\deg \nu=\deg Q=\rho(f)$, which is the required conclusion $(ii)$.\\
\textbf{Case 1.3:} If $\rho(f)>\deg Q$, then from equation \eqref{chweqsimple}
\begin{align*}
 T(r,ue^{\nu})&=T(r,f^n+wf^{n-1}f'+qe^{Q}f_c)\\
 &\leq T(r,f^{n-1})+T(r,f+wf')+T(r,q)+T(r,e^{Q})+T(r,f_c)+O(1)\\
 &\leq (n+1)T(r,f)+T(r,\frac{f+wf'}{f})+S(r,f)\\
 &=(n+1)T(r,f)+S(r,f).
\end{align*}
This implies $\deg \nu\leq \rho(f)=\deg \beta $. Now, we discuss the following sub cases:
\begin{enumerate}[(i)]
\item If $\deg \nu<\deg \beta$, we consider $\beta(z)=a_tz^{t}+a_{t-1}z^{t-1}+...+a_0, (a_t\neq 0)$ and $Q(z)=b_dz^{d}+b_{d-1}z^{d-1}+...+b_0,(b_d\neq 0)$ such that $d< \rho(f)=t$, then from equation \eqref{maineq1}, we have
$$[\alpha^{n}+w\alpha^{n-1}(\alpha'+\alpha \beta')]e^{\alpha_1}e^{na_tz^{t}}+q\alpha_ce^{\tilde{\alpha_2}}e^{a_tz^{t}}=ue^{\nu},$$
where $\tilde{\alpha_2}=Q(z)+a_{t-1}z^{t-1}+...+a_0$ is a polynomial of degree at most $t-1$. Then using Lemma \ref{imp} into the above equation, we get $u\equiv 0$, which is not possible.
\item If  $\deg \nu=\deg \beta$, then we consider $\nu(z)=\nu_{t}z^{t}+\nu_{t-1}z^{t-1}+...+\nu_{0},\nu_{t}\neq 0$. From equation \eqref{maineq1}
$$[\alpha^{n}+w\alpha^{n-1}(\alpha'+\alpha \beta')]e^{\alpha_1}e^{na_tz^{t}}+q\alpha_ce^{\tilde{\alpha_2}}e^{a_tz^{t}}=ue^{\tilde{\alpha_3}}e^{\nu_{t}z^{t}},$$
where $\tilde{\alpha_3}=\nu_{t-1}z^{t-1}+...+\nu_{0}$ is a polynomial of degree at most $t-1$.\\
If $\nu_{t}\neq na_{t}$ and $\nu_{t}\neq a_{t}$, then using Lemma \ref{imp}, $q\alpha_c\equiv 0$, which is not possible.\\
If $\nu_t=na_t$, then again using Lemma \ref{imp}, we get the same contradiction as above.\\
If $\nu_{t}=a_{t}$, then
$$[\alpha^{n}+w\alpha^{n-1}(\alpha'+\alpha \beta')]e^{\alpha_1}e^{na_tz^{t}}+(q\alpha_ce^{\tilde{\alpha_2}}-ue^{\tilde{\alpha_3}})e^{a_{t}z^{t}}\equiv 0,$$
Using Lemma \ref{imp}, $$\alpha^{n-1}(\alpha+w(\alpha'+\alpha \beta'))\equiv 0.$$
$\implies$ \qquad  \qquad \quad $\alpha+w(\alpha'+\alpha \beta')=0$.\\
This gives that $\alpha=C_3e^{-z/w-\beta}$. This implies $\deg \alpha=\deg \beta=\rho(f)$, which is not possible.
\end{enumerate} 
\textbf{Case 2:} If we consider $n\geq 2$ with $w=0$, then similarly we study 3 sub cases as for Case 1.\\
\textbf{Case 2.1:} When $\rho(f)<\deg Q$, then proceeding similar lines as in  Case 1.1, we get $\deg Q=\deg \nu$. On differentiating equation $\eqref{chweqsimple}$, we have
\begin{equation}\label{keyeq2.1}
	nf^{n-1}f'+(q'f_c+qQ'f_c+qf'_c)e^{Q}=(u'+u\nu')e^{\nu},
\end{equation}
Eliminating $e^{\nu}$ from equation \eqref{chweqsimple} and \eqref{keyeq2.1}, then we obtain
$$D_1e^{Q}+D_0=0,$$
where $$D_1=(u'+u\nu')qf_c-u(q'f_c+qQ'f_c+qf'_c),$$
$$D_0=(u'+u\nu')f^n-nuf^{n-1}f'.$$
Using Lemma \ref{imp}, we have $D_1=D_0\equiv 0$.\\
From $D_1\equiv0$, we get 
$$\frac{u'}{u}+\nu'=\frac{q'}{q}+\frac{f'_c}{f_c}+Q'.$$
On doing integration of the above equation, we have
$$qf_ce^Q=C_1ue^{\nu}.$$
where $C_1$ is a non-zero constant.\\ 
If $C_1=1$, then $qf_ce^Q=ue^{\nu}$ and by equation \eqref{chweqsimple}, we obtain
$$f^{n}\equiv 0.$$
Which is not possible. Hence conclution (i) does not hold.\\
If $C_1\neq 1$, then again proceeding same lines to Case 1.1, we get same contradiction.\\
\textbf{Case 2.2:}  When $\rho(f)=\deg Q$, Using same argument as for Case 1.2, we get required conclusion (ii).\\
\textbf{Case 2.3:} When $\rho(f)>\deg Q$, Using same argument as for Case 1.3,we get same conradiction. \\
Hence we complete the proof.
\end{proof}

\begin{proof}[\textbf{\underline{Proof of Theorem B}}]
We prove our result by contradiction. Suppose that $f$ is a finite order transcendental entire function satisfying equation \eqref{dedeg3} with a non-zero finite Borel exceptional value. Then we study following three cases:\\
\textbf{Case 1:} If $\rho(f)<1$, then using similar technique as in [\cite{chw},Theorem 1.11], we get contradiction. \\
\textbf{Case 2:}  If $\rho(f)>1$ and let $a(\neq 0)$ be a finite Borel exceptional value of $f$. Then using Weierstrass factorisation theorem, $f(z)$ can be represented in the form
\begin{equation}\label{usingwft}
f(z)=\alpha(z)e^{\beta(z)}+a,
\end{equation} 
where $\beta(z)=a_{t}z^{t}+a_{r}z^{t}+...+a_{0}; a_{t}\neq 0$ is a polynomial of degree $t>1$ and $\alpha(z)$ is an entire function such that $\rho(\alpha)<\deg \beta.$\\
From equation \eqref{dedeg3}, we have
\begin{align*}
\qquad \qquad m(r,e^{Q(z)})&=m\left(r,\frac{q(z)e^{Q(z)}f(z+c)}{q(z)f(z+c)}\right)\\
&=m\left(r,\frac{P_{1}e^{\lambda z}+P_{2}e^{-\lambda z}-f^{2}(z)(f(z)+wf^{'}(z))}{q(z)f(z+c)}\right)\\
&\leq m\left(r,\frac{P_{1}e^{\lambda z}+P_{2}e^{-\lambda z}}{q(z)f(z+c)}\right)+ m\left(r,\frac{f^{2}(z)(f(z)+wf^{'}(z))}{q(z)f(z+c)}\right)+O(1)\\
&\leq  m\left(r,\frac{1}{q(z)f(z+c)}\right)+m(r,P_{1}e^{\lambda z}+P_{2}e^{-\lambda z})+m(r,f^{2})+\\
&m\left(r,\frac{1}{q(z)}\right)+m\left(r,\frac{f(z)}{q(z)f(z+c)}\right)+m\left(r,\frac{f'(z)}{q(z)f(z+c)}\right)+O(1)
\end{align*}
Applying Lemma \ref{il} and \ref{hk} into the above inequality, we have
\begin{align*}
T(r,e^{Q(z)})&\leq  m\left(r,\frac{1}{q(z)f(z+c)}\right)+2T(r,f)+m\left(r,\frac{f'(z)}{q(z)f(z+c)}\right)+S(r,f)\\
&\leq m\left(r,\frac{f(z)}{q(z)f(z+c)}\right)+m\left(r,\frac{1}{f(z)}\right)+2T(r,f)+m\left(r,\frac{f'(z)f(z)}{q(z)f(z)f(z+c)}\right)\\
& \qquad \qquad  \qquad \qquad \qquad +S(r,f)\\
&\leq 3T(r,f)+S(r,f).
\end{align*}
This implies that $\deg Q\leq \rho(f)=t$ and we also know that $\deg Q\geq 1$.\\
First suppose $1\leq \deg Q<t$, then substituting the value of $f(z)$ from equation \eqref{usingwft} into equation \eqref{dedeg3}, we have
\begin{equation}\label{keyeq3.1}
\begin{split}
A_{1}(z)e^{3\beta(z)}+A_{2}(z)e^{2\beta(z)}&+A_{3}(z)e^{\beta(z)}+aq(z)e^{Q(z)}+\\
&q(z)\alpha(z+c)e^{Q(z)+\beta(z+c)}-P_{1}e^{\lambda z}-P_{2}e^{-\lambda z}+a^{3}=0,
\end{split}
\end{equation}

where
\begin{align*}
A_{1}(z)&=\alpha^{3}(z)+w\alpha^{3}(z)\beta^{'}(z)+w\alpha^{2}(z)\alpha^{'}(z)\\
A_{2}(z)&=3a\alpha^{2}(z)+2wa\beta^{'}(z)\alpha^{2}(z)+2wa\alpha(z)\alpha^{'}(z)\\
A_{3}(z)&=3a^{2}\alpha(z)+wa^{2}\alpha(z)\beta^{'}(z)+wa^{2}\alpha^{'}(z).
\end{align*}
On simplifying equation \eqref{keyeq3.1}, we get
\begin{equation}\label{keyeq3.1.1}
\begin{split}
A_{1}(z)e^{3\beta(z)}+A_{2}(z)e^{2\beta(z)}&+(A_{3}(z)e^{\beta_{1}(z)}+q(z)\alpha(z+c)e^{\delta(z)})e^{a_{t}z^{t}}+\\
	&(aq(z)e^{Q(z)}-P_{1}e^{\lambda z}-P_{2}e^{-\lambda z}+a^{3})=0,
\end{split}
\end{equation}
where $\beta_{1}(z)=a_{t-1}z^{t-1}+a_{t-2}z^{t-2}+...+a_{0}$ and $\delta(z)=Q(z)+a_{t}[\binom{t}{1}cz^{t-1}+\binom{t}{2}c^{2}z^{t-2}+...+c^{t}]+a_{t-1}(z+c)^{t-1}+...+a_{0}$ are polynomials of degree at most $t-1$. Applying Lemma \ref{imp} into equation \eqref{keyeq3.1.1}, we have 
\begin{equation}\label{subkeyeq3.1.1}
aq(z)e^{Q(z)}-P_{1}e^{\lambda z}-P_{2}e^{-\lambda z}+a^{3}\equiv 0
\end{equation}
Again applying Lemma \ref{imp} into equation \eqref{subkeyeq3.1.1}, we have the following sub cases:
\begin{enumerate}[(i)]
\item If $Q(z)\neq \pm \lambda z$, then we have $P_{1}\equiv 0\equiv P_{2}=a^{3}$, which is not possible.
\item  If $Q(z)=\lambda z$, then we have $P_{2}\equiv 0\equiv a^{3}$, which is not possible
\item If $Q(z)=-\lambda z$, then we have $P_{1}\equiv 0\equiv a^{3}$, which is not possible.
\end{enumerate}
Now suppose $\deg Q=t$, and let $Q(z)=b_{t}z^{t}+b_{t-1}z^{t-1}+...+b_{0}; b_{t}\neq 0$ is a polynomial. Then from equation \eqref{keyeq3.1}, we have
\begin{equation}\label{subkeyeq3.1.2}
\begin{split}
A_{1}(z)e^{3\beta_{1}(z)}e^{3a_{t}z^{t}}&+A_{2}(z)e^{2\beta_{1}(z)}e^{2a_{t}z^{t}}+A_{3}(z)e^{\beta_{1}(z)}e^{a_{t}z^{t}}+aq(z)e^{Q_{1}(z)}e^{b_{t}z^{t}}+\\
&q(z)\alpha(z+c)e^{Q_{1}(z)+\beta_{1}(z+c)}e^{b_{t}z^{t}+a_{t}(z+c)^{t}}-(P_{1}e^{\lambda z}+P_{2}e^{-\lambda z}-a^{3})=0,
\end{split}
\end{equation}
where $\beta_{1}(z)=a_{t-1}z^{t-1}+a_{t-2}z^{t-2}+...+a_{0}$ and $Q_{1}(z)=b_{t-1}z^{t-1}+b_{t-2}z^{t-2}+...+b_{0}$ are polynomials of degree at most $t-1$. Next we discuss the following sub cases:
\begin{enumerate}[(i)]
\item If $b_{t}\neq ka_{t}$ for any $k=1,2,3$, and $b_{t}\neq -a_{t}$ then applying Lemma \ref{imp}, we have $aq(z)e^{Q_{1}(z)}\equiv 0$, this implies either $a\equiv 0$ or $q(z)\equiv 0$, which is not possible.
\item If $b_{t}=ka_{t}$ for some $k=1,2,3$, say $b_{t}=a_{t}$, then equation \eqref{subkeyeq3.1.2} becomes
\begin{equation*}
\begin{split}
A_{1}(z)e^{3\beta_{1}(z)}e^{3a_{t}z^{t}}&+(A_{2}(z)e^{2\beta_{1}(z)}+q(z)\alpha(z+c)e^{\tilde{\delta}(z)})e^{2a_{t}z^{t}}+\\
&(A_{3}(z)e^{\beta_{1}(z)}+aq(z)e^{Q_{1}(z)})e^{a_{t}z^{t}}-(P_{1}e^{\lambda z}+P_{2}e^{-\lambda z}-a^{3})=0,
\end{split}
\end{equation*}
where $\tilde{\delta}(z)=Q_{1}(z)+\beta_{1}(z+c)+a_{t}[\binom{t}{1}cz^{t-1}+\binom{t}{2}c^{2}z^{t-2}+...+c^{t}]$ is a polynomial of degree at most $t-1$. Applying Lemma \ref{imp} into the above equation, we have
$P_{1}e^{\lambda z}+P_{2}e^{-\lambda z}-a^{3}\equiv 0$, which is not possible.
\item If $b_{t}=-a_{t}$, then equation \eqref{subkeyeq3.1.2} becomes
\begin{equation*}
\begin{split}
A_{1}(z)e^{3\beta_{1}(z)}e^{3a_{t}z^{t}}&+A_{2}(z)e^{2\beta_{1}(z)}e^{2a_{t}z^{t}}+A_{3}(z)e^{\beta_{1}(z)}e^{a_{t}z^{t}}+aq(z)e^{Q_{1}(z)}e^{b_{t}z^{t}}+\\
&(q(z)\alpha(z+c)e^{\tilde{\delta}(z)}-P_{1}e^{\lambda z}+P_{2}e^{-\lambda z}-a^{3})=0,
\end{split}
\end{equation*}
Applying Lemma \ref{imp} into the above equation, we have $aq(z)e^{Q_{1}(z)}\equiv 0$, this implies either $a\equiv 0$ or $q(z)\equiv 0$, which is not possible.
\end{enumerate}
\textbf{Case 3:} If $\rho(f)=1$, then we can represent $f(z)$ as 
\begin{equation}
f(z)=\gamma(z)e^{bz+d}+a,
\end{equation}
where $b\neq 0, d\in \mathbb{C}$ and $\gamma(z)$ is an entire such that $\rho(\gamma)<1$.\\
Now proceeding to similar lines as in case of $\rho(f)>1$, we have $\deg Q=1$. Thus we can put $Q(z)=pz+s$,  $(0\neq p,s\in\mathbb{C})$ in the equation \eqref{dedeg3}, we have
\begin{equation}\label{keyeq3.3}
B_{1}(z)e^{3bz}+B_{2}(z)e^{2bz}+B_{3}(z)e^{bz}+B_{4}(z)e^{(p+b)z}+aq(z)e^{s}e^{pz}+a^{3}-P_{1}e^{\lambda z}-P_{2}e^{-\lambda z}=0,
\end{equation}
where
\begin{align*}
B_{1}(z)&=(\gamma^{3}(z)+w\gamma^{2}(z)\gamma^{'}(z)+wb\gamma^{3}(z))e^{3d}\\
B_{2}(z)&=(3a\gamma^{2}(z)+2aw\gamma(z)\gamma^{'}(z)+2abw\gamma^{2}(z))e^{2d}\\
B_{3}(z)&=(3a^{2}\gamma(z)+a^{2}w\gamma^{'}(z)+a^{2}bw\gamma(z))e^{d}\\
B_{4}(z)&=q(z)\gamma(z+c)e^{bc+s+d}.
\end{align*}
Here we discuss the following sub cases:
\begin{enumerate}[(i)]
\item If $kb\neq p$ and $kb\neq \pm\lambda$ for any $k=1,2,3$, then applying Lemma \ref{imp} to equation \eqref{keyeq3.3}, we have $P_{1}\equiv 0\equiv P_{2}$, which is not possible.
\item If $kb=p$ for some $k=1,2,3$ and $kb\neq \pm\lambda$ for any $k=1,2,3$, then applying Lemma \ref{imp} to equation \eqref{keyeq3.3}, we obtain same contradiction as for previous case.
\item If  $kb\neq p$ for any $k=1,2,3$ and $kb=\lambda$ for some $k=1,2,3$, then applying Lemma \ref{imp} to equation \eqref{keyeq3.3}, we have $P_{2}\equiv 0$, which is not possible.
\item If  $kb\neq p$ for any $k=1,2,3$ and $kb=-\lambda$ for some $k=1,2,3$, then applying Lemma \ref{imp} to equation \eqref{keyeq3.3}, we have $P_{1}\equiv 0$, which is not possible.
\item  If $kb=p$ and $kb=\lambda$ for some $k=1,2,3$, then applying Lemma \ref{imp} to equation \eqref{keyeq3.3}, we have $P_{2}\equiv 0\equiv a^{3}$, which is not possible.
\item If $kb=p$ and $kb=-\lambda$ for some $k=1,2,3$, then again applying Lemma \ref{imp} to equation \eqref{keyeq3.3}, we have $P_{1}\equiv 0\equiv a^{3}$, which is not possible.
\end{enumerate}
Thus we complete the proof.\\
\end{proof}

\begin{proof}[\textbf{\underline{Proof of Theorem C}}]
Suppose that $f$ is a transcendental meromorphic function satisfying equation \eqref{dedeg2} and given that $f$ has finitely many poles, then from equation \eqref{dedeg2}, and Lemma \ref{il} and \ref{whllem},
\begin{align*}
nT(r,f)=T(r,f^{n})&=T(r,P_{1}e^{\alpha_{1}(z)}+P_{2}e^{\alpha_{2}(z)}+P_{3}e^{\alpha_{3}(z)}-P(z,f))\\
&\leq T(r,P_{1}e^{\alpha_{1}(z)}+P_{2}e^{\alpha_{2}(z)}+P_{3}e^{\alpha_{3}(z)})+T(r,P(z,f))+\log 2\\
&\leq T(r,P_{1}e^{\alpha_{1}(z)}+P_{2}e^{\alpha_{2}(z)}+P_{3}e^{\alpha_{3}(z)})+m\left(r,\frac{P(z,f)}{f}f\right)+\\
&\qquad  \qquad \qquad \qquad N(r,P(z,f)) +\log 2\\
&\leq Ar^{m}+T(r,f)+S(r,f),
\end{align*}
where $$A=\frac{\text{sum of the leading coefficients of} \alpha_{1}(z), \alpha_{2}(z)\text{and}\alpha_{3}(z)}{\pi}$$and $$m=\max\{\deg \alpha_{1}(z),\deg \alpha_{2}(z),\deg \alpha_{3}(z)\}.$$

Thus we have $(n-1)T(r,f)\leq Ar^{m}+S(r,f)$ which implies that $f$ has finite order.\\
Now we prove rest part of the theorem by contradiction. Suppose $f(z)=s(z)e^{q(z)}$ is a solution of equation \eqref{dedeg2}, where $s(z)$ is a non-zero rational function and $q(z)$ is a non constant polynomial, then from equation \eqref{dedeg2}, we get
\begin{align*}
P_{1}(z)e^{\alpha_{1}(z)}+P_{2}(z)e^{\alpha_{2}(z)}+P_{3}(z)e^{\alpha_{3}(z)}&=f^{n}(z)+P(z,f)\\&=f^{n}(z)+\sum_{i=0}^{i=l}a_{i}(z)f^{(i)}(z)+\tilde{a_{0}}(z)\\
&=(s(z)e^{q(z)})^{n}+\sum_{i=0}^{i=l}a_{i}(z)(s(z)e^{q(z)})^{(i)}+\tilde{a_{0}}(z)
\end{align*} 
This implies 
\begin{equation}\label{keyeqlast}
P_{1}(z)e^{\alpha_{1}(z)}+P_{2}(z)e^{\alpha_{2}(z)}+P_{3}(z)e^{\alpha_{3}(z)}-s^{n}(z)e^{nq(z)}-Q(z)e^{q(z)}-\tilde{a_{0}}(z)=0,
\end{equation}
where $Q(z)=a_{0}(z)s(z)+a_{1}(z)[s^{'}(z)+s(z)q^{'}(z)]+...+a_{l}(z)[s^{l}(z)+...+s(z)(q^{'}(z))^{l}]$ is a rational function. Now we discuss the following cases: 
\begin{enumerate}[(i)]
\item If $q(z)-\alpha_{j}(z)\neq $  constant, and $nq(z)-\alpha_{j}(z)\neq $ constant, for any $j=1,2,3.$ Then applying Lemma \ref{imp} into equation \eqref{keyeqlast}, we have $P_{1}(z)\equiv 0\equiv P_{2}(z)=P_{3}(z)$, which is not possible. 
\item If $q(z)-\alpha_{j}(z)=$ constant, for some $j=1,2,3$, say $q(z)-\alpha_{1}(z)=$ constant, and $nq(z)-\alpha_{j}(z)\neq$ constant, for any $j=1,2,3$. Let $q(z)=b_{r}z^{r}+b_{r-1}z^{r-1}+...+b_{0}$ and $\alpha_{1}(z)=b_{r}z^{r}+b_{r-1}z^{r-1}+...+c_{0}$. Then $q(z)-\alpha_{1}(z)=b_{0}-c_{0}$ and equation \eqref{keyeqlast} becomes
$$(P_{1}(z)e^{c_{0}}-Q(z)e^{b_{0}})e^{\tilde{\alpha_{1}}(z)}+P_{2}(z)e^{\alpha_{2}(z)}+P_{3}(z)e^{\alpha_{3}(z)}-s^{n}(z)e^{nq(z)}-\tilde{a_{0}}(z)=0,$$where $\tilde{\alpha_{1}}(z)=b_{r}z^{r}+b_{r-1}z^{r-1}+...+b_{1}z.$
Applying Lemma \ref{imp} to the above equation, we have $P_{2}(z)\equiv 0\equiv P_{3}(z)$, which is not possible.
\item If $q(z)- \alpha_{j}(z)\neq$ constant, for any $j=1,2,3$ and $nq(z)- \alpha_{j}(z)=$ constant, for some $j=1,2,3$, say $nq(z)-\alpha_{1}(z)=$ constant.  Then applying Lemma \ref{imp} into equation \eqref{keyeqlast}, we get the same contradiction as for just previous case.
\item If $q(z)-\alpha_{j}(z)=$ constant, for some $j=1,2,3$, say $q(z)-\alpha_{1}(z)=$ constant and $nq(z)-\alpha_{k}(z)=$ constant, for some $k$, say $nq(z)-\alpha_{2}(z)=$ constant. Then applying Lemma \ref{imp} into equation \eqref{keyeqlast}, we have $P_{3}(z)\equiv 0$, which is not possible.
\end{enumerate}
Hence we complete the proof.
\end{proof}

\begin{proof}[\textbf{\underline{Proof of Theorem D}}]
Since $f$ is an entire function, then from equation \eqref{dedeg3sp} and Lemma \ref{il}, we have
 \begin{align*}
 T(r,P_{1}e^{\alpha_{1}z}+P_{2}e^{\alpha_{2}z}+P_{3}e^{\alpha_{3}z})&=T(r,f^{n}+Q_{2}^{*}(z,f))\\
  &\leq nT(r,f)+T\left(r,\frac{Q_{2}^{*}(z,f)}{f^{2}}f^{2}\right)\\
 &\leq nT(r,f)+T\left(r,\frac{Q_{2}^{*}(z,f)}{f^{2}}\right)+2T(r,f)\\
 &\leq (n+3)T(r,f)+S(r,f)
 \end{align*}
As we know that $\rho(P_{1}e^{\alpha_{1}z}+P_{2}e^{\alpha_{2}z}+P_{3}e^{\alpha_{3}z})=1$ by using Lemma \ref{whllem}, thus we get $\rho(f)\geq 1$.\\
First suppose $\rho(f)>1$ and given that $f$ is an entire satisfying equation \eqref{dedeg3sp} with $\lambda(f)<\rho(f)$. Then using Weierstrass factorization theorem, $f$ must be a transcendental entire function in the form of
\begin{equation*}\label{wfteqq}
	f(z)=H(z)e^{\beta(z)},
\end{equation*} 
where $\beta(z)$ is a polynomial of degree $r>1$ and $H(z)$ is an entire function such that $\rho(H)<r$. Now from equation \eqref{dedeg3sp}, we have
\begin{align*}
H^{n}(z)e^{n\beta(z)}+\sum_{i=1}^{i=2}b_{i}(z)e^{i\beta(z)}+a_{0}(z)=P_{1}e^{\alpha_{1}z}+P_{2}e^{\alpha_{2}z}+P_{3}e^{\alpha_{3}z},
\end{align*}
where $b_{1}(z),b_{2}(z)$ are meromorphic functions of order $<r$ and $a_{0}(z)$ is a rational function. This implies
\begin{equation*}
H^{n}(z)e^{n\beta(z)}+b_{2}(z)e^{2\beta(z)}+b_{1}(z)e^{\beta(z)}+(a_{0}(z)-P_{1}e^{\alpha_{1}z}-P_{2}e^{\alpha_{2}z}-P_{3}e^{\alpha_{3}z})=0
\end{equation*} 
Applying Lemma \ref{imp} into the above equation, we have
\begin{equation}\label{keyeqD1}
a_{0}(z)-P_{1}e^{\alpha_{1}z}-P_{2}e^{\alpha_{2}z}-P_{3}e^{\alpha_{3}z}\equiv 0,
\end{equation}
which is not possible. Otherwise, again applying Lemma \ref{imp} into equation \eqref{keyeqD1}, we get $P_{1}\equiv 0\equiv P_{2}=P_{3}$, which is contradiction. Thus $\rho(f)=1$.\\
Now considering 
\begin{equation}\label{eqD2}
f(z)=H(z)e^{az+b},
\end{equation}
where  $0\neq a,b\in\mathbb{C}$ and $H(z)$ is an entire function of order $<1$. From equation \eqref{dedeg3sp}, we have
\begin{equation*}
H^{n}(z)e^{n(az+b)}+b_{2}(z)e^{2(az+b)}+b_{1}(z)e^{az+b}+a_{0}(z)-P_{1}e^{\alpha_{1}z}-P_{2}e^{\alpha_{2}z}-P_{3}e^{\alpha_{3}z}=0
\end{equation*}
Applying Lemma \ref{imp} into the above equation, there must exist positive integers $k_{1},k_{2}$ and $k_{3}$ with $\{k_{1},k_{2},k_{3}\}=\{1,2,3\}$ such that
$na=\alpha_{k_{1}}, 2a=\alpha_{k_{2}}, a=\alpha_{k_{3}}$ which gives $\alpha_{k_{1}}:\alpha_{k_{2}}:\alpha_{k_{3}}=n:2:1$. We also have 
$H^{n}(z)e^{nb}=P_{k_{1}}$, this gives $H(z)\equiv P_{k_{1}}^{1/n}e^{-b}$ is a constant, $b_{2}(z)e^{2b}\equiv P_{k_{2}}$, $b_{1}(z)e^{b}\equiv P_{k_{3}}$, and $a_{0}(z)\equiv 0$. Substituting the value of $H(z)$ into equation \eqref{eqD2}, we have $f(z)=P_{k_{1}}^{1/n}e^{az}=Ce^{az}$, where $C$ is a non-zero constant.
\end{proof}


\begin{thebibliography}{99}
\bibitem{cl} Chen, J. F., Lian, G., \textit{Expressions of meromorphic solutions of a certain type of nonlinear complex differential equations}, Bulletin of the Korean Mathematical Society, 57(4) (2020), 1061-1073.
\bibitem{ilaine} Laine, I., \textit{Nevanlinna Theory and Complex Differential Equations}, W. de Gruyter, Berlin (1993).
\bibitem{hayman} Hayman, W.K.,\textit{Meromorphic Functions}, Clarendon Press, Oxford (1964)
\bibitem{hk1} Halburd, R.G. and Korhonen, R.J., \textit{Difference analogue of the lemma on the logarithmic derivative with applications to difference equations}, J. Math. Anal. Appl., 314 pp.477–487 (2006).
\bibitem{hk2} Halburd, R.G. and Korhonen, R.J., Nevanlinna theory for the difference operator, \textit{Ann. Acad. Sci. Fenn. Math.}, \textbf{31} (2006), 463–478.
\bibitem{hkt} Halburd, R., Korhonen, R. and Tohge, K., Holomorphic curves with shift-invariant hyperplane preimages, \textit{Transactions of the American Mathematical Society}, \textbf{366}, no. 8 (2014), 4267-4298.
\bibitem{xue} Xue, B., Entire solutions of certain type of non-linear differential equations. Mathematica Slovaca, 70(1) (2020), 87-94.
\bibitem{mirsky} L. Mirsky, An Introduction to Linear Algebra, Oxford, at the Clarendon Press, 1955.
\bibitem{li} Li, P., \textit{Entire solutions of certain type of differential equations}, Journal of mathematical analysis and applications, 344(1), pp.253-259 (2008).
\bibitem{li1} Li P., \textit{Entire solutions of certain type of differential equations II}, Journal of mathematical analysis and applications, 375(1), pp.310-319 (2011).
\bibitem{ly} Li, P. and Yang, C.C., On the nonexistence of entire solutions of a certain type of nonliear differential equations, \textit{J. Math. Anal. Appl.}, \textbf{320}, 827-835 (2006).
\bibitem{yybook} C.C.Yang, H.X.Yi: \textit{Uniqueness Theory of Meromorphic Functions}, Science Press, Kluwer Academic, Dordrecht, Beijing, (2003).
\bibitem{chw}Chen, W., Hu, P. and Wang, Q.: \textit{Entire Solutions of Two Certain Types of Non-linear Differential-Difference Equations}. Computational Methods and Function Theory, 21(2), pp.199-218 (2021).
\bibitem{cf} Chiang, Y.M. and Feng, S.J.: \textit{On the Nevanlinna characteristic of f(z+$\eta$) and difference equations in the complex plane}, The Ramanujan Journal, 16(1), pp.105-129 (2008).
\bibitem{rx} Rong J, Xu J., \textit{Three results on the nonlinear differential equations and differential-difference equations}, Mathematics, 539 (2019).
\bibitem{ww} Wang Q, Wang Q.: \textit{Study on the existence of solutions to two specific types of differential-difference equations}, Turkish Journal of Mathematics, 43(2), pp.941-54 (2019).
\bibitem{whl} Wen, Z.T., Heittokangas, J, Lain, I, \textit{Exponential polynomials as solutions of certain nonlinear difference equations}. Acta Mathematica Sinica, English Series, 28(7),pp.1295-306 (2012).
\bibitem{yanglo} Yang Lo, \textit{Value Distribution Theory}, Translated and revised from the 1982 Chinese Original, Springer-Verlag, Berlin (1993).
\end{thebibliography}
\end{document}